\documentclass[12pt]{article}

\usepackage{amsmath}
\usepackage{amsfonts}
\usepackage{amssymb}
\usepackage{graphicx}

\newtheorem{proposition}{Proposition}[section]
\newtheorem{theorem}{Theorem}[section]
\newtheorem{lemma}[theorem]{Lemma}
\newtheorem{corollary}[theorem]{Corollary}
\newtheorem{remark}[theorem]{Remark}

\def\phi{{\varphi}}

\DeclareSymbolFont{AMSb}{U}{msb}{m}{n}
\DeclareMathSymbol{\N}{\mathbin}{AMSb}{"4E}
\DeclareMathSymbol{\Z}{\mathbin}{AMSb}{"5A}
\DeclareMathSymbol{\R}{\mathbin}{AMSb}{"52}
\DeclareMathSymbol{\Q}{\mathbin}{AMSb}{"51}
\DeclareMathSymbol{\I}{\mathbin}{AMSb}{"49}
\DeclareMathSymbol{\C}{\mathbin}{AMSb}{"43}

\begin{document}

\title{A convergent  algorithm for the hybrid problem of reconstructing conductivity from minimal interior data }

\author{{Amir Moradifam\footnote{Department of Mathematics, University of Toronto, Toronto, Ontario, Canada M5S 2E4. E-mail: amir@math.toronto.edu. The author is supported by a MITACS Postdoctoral Fellowship. }
\qquad Adrian Nachman\footnote{Department of Mathematics and the
Edward S. Rogers Sr. Department of Electrical and Computer
Engineering, University of Toronto, Toronto, Ontario, Canada M5S
2E4. E-mail: nachman@math.toronto.edu. The author is supported in part by an NSERC Discovery Grant.}\qquad Alexandre Timonov \footnote{Division of Mathematics and Computer Science, University of South Carolina Upstate, Spartanburg, SC, USA.  E-mail: atimonov@uscupstate.edu.}}}
\date{\today}

\smallbreak \maketitle

\begin{abstract}
We consider the hybrid problem of reconstructing the isotropic electric conductivity of a body $\Omega$ from interior Current Density Imaging data obtainable using MRI measurements. We only require knowledge of the magnitude $|J|$ of one current generated by a given voltage $f$ on the boundary  $\partial\Omega$. As previously shown, the corresponding voltage potential u in $\Omega$ is a minimizer of the weighted least gradient problem

\[u=\hbox{argmin} \{\int_{\Omega}a(x)|\nabla u|: u \in H^{1}(\Omega),
\ \ u|_{\partial \Omega}=f\},\]
with $a(x)= |J(x)|$. In this paper we present an alternating split Bregman algorithm for treating such least gradient problems, for $a\in L^2(\Omega)$ non-negative and  $f\in H^{1/2}(\partial \Omega)$. 
We give a detailed convergence proof by focusing to a large extent on the dual problem. This leads naturally to the alternating split Bregman algorithm. The dual problem also turns out to yield a novel method to recover the full vector field $J$ from knowledge of its magnitude, and of the voltage $f$ on the boundary. We then present several numerical experiments that illustrate the convergence behavior of the proposed algorithm.

\end{abstract}
\maketitle

\section{Introduction}

The classical Electrical Impedance Tomography (EIT) problem seeks to obtain quantitative information on the electrical conductivity $\sigma$ of a body from multiple measurements of voltages and corresponding currents at its surface. The extensive study of this problem has led to major mathematical advances on uniqueness and reconstruction methods for Inverse Problems with boundary data. See the excellent reviews \cite{isaacsonReview}, \cite{borcea}, and \cite{GLKU}.
However, by now it is well understood that the problem is exponentially ill-posed, yielding images of low resolution away from the boundary \cite{isacson}, \cite{Man}.

A new class of Inverse Problems seeks to significantly improve both the
quantitative accuracy and the resolution of traditional Inverse Boundary Value Problems by using data which can be
determined in the interior of the object. These have been dubbed Hybrid Methods, as they usually involve the combination of two different kinds of physical measurements, and several recent advances are presented in the present issue of the journal.

In this paper we continue our study of the Current Density Impedance Imaging (CDII) problem of reconstructing the conductivity of a body based on  measurement of currents in its interior. Such measurements have been possible since the early 1990 due to the pioneering work of M. Joy's group at the University of Toronto \cite{joy89}, \cite{joy}. The idea was to use Magnetic Resonance Imaging (MRI) in a novel way, to determine the magnetic flux density $B$ induced by an applied current. It is important to note that in the problem addressed in this paper, we only require knowledge of the magnitude $|J|$ of one current generated by a given voltage $f$ on the boundary  $\partial\Omega$. The analytic and numerical methods presented here do not necessarily depend on MRI. Since the results only require knowledge of the magnitude of one current, they may lead to simpler physical methodologies to obtain such data (see for instance \cite{Yuan07}). 

The problem of recovering the isotropic conductivity $\sigma$ of an object from knowledge of the magnitude of one current density $|J|$ in the interior was studied in \cite{NTT07, NTT08, NTT10, MNT}. See \cite{NTT11} for a review, and for numerous references to other hybrid approaches to conductivity imaging. In this paper we present an alternating split Bregman algorithm for the numerical solution of this problem, along with a convergence proof.

Let $\sigma$ be the the isotropic conductivity of an object $\Omega\subset \R^n$, $n\geq 2$ and let $J$ be the current density vector field generated by imposing a given boundary voltage $f$. Then the corresponding voltage potential $v$ satisfies the elliptic equation
\begin{equation}
\nabla \cdot \left(\sigma \nabla v \right)=0, \ \ v|_{\partial \Omega}
=f.\end{equation}
By Ohm's law $ J=-\sigma \nabla v$. Hence the voltage potential $v$ satisfies the degenerate elliptic equation 
\begin{equation}\label{Deg}
\nabla \cdot \left( \frac{|J|}{|\nabla v|}\nabla v\right)=0, \ \ v|_{\partial \Omega}=f. 
\end{equation}
 In general there is no uniqueness for viscosity solutions \cite{CIL} of the above equation \cite{SZ, NTT08}. However, in  \cite{NTT08} and \cite {MNT} authors proved that if the potential $v$ satisfies  eqauation (1) above then it minimizes the energy functional $E(v)=\int_{\Omega} |J||\nabla v|$ associated to the equation (\ref{Deg}) and, moreover that this functional has a unique minimizer in $H^1(\Omega)$.  Thus, the conductivity is uniquely determined by the magnitude of the current generated by imposing a given boundary voltage and the corresponding voltage potential is the unique solution of the (infinite-dimensional) minimization problem 
\begin{equation}\label{min_prob}
\hbox{min} \{\int_{\Omega}|J||\nabla v|: v \in H^{1}(\Omega).
\ \ v|_{\partial \Omega}=f\}.
\end{equation}
A simple iterative procedure to solve (\ref {min_prob}) was given in \cite{NTT08}. This method was only defined when Dirichlet problems such as (1) for successive approximations of $\sigma$ were guaranteed to produce solutions with non-vanishing gradients in $\Omega$. For planar regions, this led to the requirement that the given boundary voltage
be almost two-to-one. (See \cite{NTT08} for the precise definition).
 
In this paper we present a convergent alternating split Bregman algorithm to find the minimizer of (\ref{min_prob}) for any given boundary voltage $f\in H^{1/2}(\partial \Omega)$.
More generally, let $\Omega$ be a bounded region in $\R^n$, $n \geq 2$.  Also let $a\in L^2(\Omega)$ be a non-negative function and let $f\in H^{1/2}(\partial\Omega)$. Consider the minimization problem 
\begin{equation}\label{GMinProb}
\min_{\{v\in H^1(\Omega),  v|_{\partial \Omega}=f\}} \int_{\Omega}a(x)|\nabla v|,
\end{equation}
and assume that it has an optimal solution in $H^1(\Omega)$. The algorithm we present in this paper will converge to a minimizer of(\ref{GMinProb}).  

The problem (\ref{GMinProb}) belongs to a general class of problems of the form 
\begin{equation}\label{oldP}
\hspace{1cm}\min_{u \in H_1}G(u)+ F(L u),
\end{equation}
where $L: H_1\rightarrow H_2$ is a bounded linear operator, the functions $G:H_1\rightarrow \R \cup \{\infty \}$ and $F:H_2\rightarrow \R \cup \{\infty\}$ are proper, convex and lower semi-continuous, and $H_1$ and $H_2$ are real Hilbert spaces. To see this, let $H_1=H^1_0(\Omega)$, $H_2=(L^2(\Omega))^n$, and $Lu=\nabla u$; fix $u_f\in H^1(\Omega)$ with $u_f|_{\Omega}=f$ and define $F:(L^2(\Omega))^n\rightarrow \R$ and $G: H^1_0(\Omega)\rightarrow \R$ as follows
\begin{equation}\label{FandG}
F(d):=\int_{\Omega} a|d+\nabla u_f|dx, \ \ G(u)\equiv 0.
\end{equation}

One approach to the problem (\ref{oldP}) is to write it as a constrained minimization problem
\begin{equation}
\min_{u\in H_1, d\in H_2}  G(u)+F(d) \ \ \ \ \hbox{subject to} \ \ \ \ Lu=d,
\end{equation}
which leads to the unconstrained problem:
\begin{equation}
\min_{u\in H_1, d\in H_2}  G(u)+F(d) +\frac{\lambda}{2} \|Lu-d\|^2. 
\end{equation}
To solve the above problem, Goldstein and Osher \cite{GO} introduced the split Bregman method: 
\begin{eqnarray}\label{SB}
(u^{k+1},d^{k+1})&=&\hbox{argmin}_{u\in H_1, d\in H_2} \{G(u)+F(d)+\frac{\lambda}{2}\parallel b^k+Lu -d\parallel^2_{2}\}, \\
b^{k+1}&=& b^k+Lu^{k+1}-d^{k+1}.\nonumber
\end{eqnarray}
Since the joint minimization problem (\ref{SB}) in both $u$ and $d$ could sometimes be hard or expensive to solve exactly, Goldstein and Osher \cite{GO} proposed the following 
\textit{alternating split Bregman algorithm} to solve the problem $(\ref{oldP})$
\begin{eqnarray}
u^{k+1}&=&\hbox{argmin}_{u \in H_1} \{G(u)+\frac{\lambda}{2}\parallel b^k+Lu^{k+1}-d\parallel^{2}_2\},\\
d^{k+1}&=&\hbox{argmin}_{d \in H_2} \{F(d)+\frac{\lambda}{2}\parallel b^k+Lu^{k+1}-d\parallel^{2}_2\},\\
b^{k+1}&=&b^k+Lu^{k+1}-d^{k+1}.
\end{eqnarray}

Cai, Osher, and Shen \cite{COS} and independently Setzer \cite{Setzer, Setzer1} proved convergence results for the above algorithm, under the assumption that  $H_1$ and $H_2$ are finite dimensional.  Motivated by the infinite dimensional problem (\ref{GMinProb}), in \cite{MN} the authors recently presented general convergence results for the alternating split Bregman algorithm in infinite dimensional Hilbert spaces. In this paper we will study the following alternating split Bregman algorithm for the Dirichlet problem (\ref{GMinProb}).  \\

 {\bf Algorithm 1 (Alternating split Bregman algorithm for weighted least gradient Dirichlet problems)}\\

Let $u_f \in H^1(\Omega)$ with $u_f|_{\partial \Omega}=f$  and initialize $b^0, d^0 \in (L^2(\Omega))^n$. For $k\geq 1$:\\ 
\begin{enumerate}
\item Solve \[ \Delta u^{k+1}= \nabla \cdot (d^k(x)-b^k(x)), \ \ u^{k+1}|_{\partial
\Omega}=0.\] 
\item Compute
 \small{
\[d^{k+1}:=\left\{ \begin{array}{ll}
\max \{|\nabla u^{k+1}+\nabla u_f +b^k|-\frac{a}{\lambda},0\}\frac{\nabla u^{k+1}+\nabla u_f +b^k}{|\nabla u^{k+1}+\nabla u_f  +b^k|} -\nabla u_f \ \ &\hbox{if} \ \ |\nabla u^{k+1}(x)+ \nabla u_f +b^k(x)|\neq 0,\\
-\nabla u_f   \ \ &\hbox{if} \ \ |\nabla u^{k+1}(x)+ \nabla u_f +b^k(x)|=0.
\end{array} \right.
\] 
 }
\item Let
\[b^{k+1}(x)=b^{k}(x)+\nabla u^{k+1}(x)-d^{k+1}(x).\]\\
\end{enumerate}

The following theorem, which we will prove in Section 3, guarantees convergence of the Algorithm 1.

\begin{theorem}\label{theoMain}
Let $a \in L^2(\Omega)$ be a non-negative function and  $f\in H^{1/2}(\partial \Omega)$. Then for any $u_f \in H^1(\Omega)$ and any $b^0, d^0 \in (L^2(\Omega))^n$ the sequences $\{b^k\}_{k\in N}$, $\{d^k\}_{k\in N}$, and $\{u^k\}_{k\in N}$ produced by Algorithm 1 converge weakly to some $\hat{b}$, $\hat{d}$, and $\hat{u}$. Moreover  $\{\nabla u^{k+1}-d^k\}_{k\in N}$ converges strongly to zero and
\begin{equation}\label{conv-est}
\sum_{k=0}^{\infty}\parallel\nabla u^{k+1}-d^k \parallel ^2<\infty.
\end{equation}
Furthermore $\bar{u}:=\hat{u}+u_f$ is a solution of the minimization problem ($\ref{GMinProb}$),  $\hat{d}=\nabla \hat{u}$, $\nabla \cdot \hat{b} \equiv 0$, and

\begin{equation}\label{bhat}
\hat{b}=\frac{1}{\lambda} \left( a\frac{\nabla \bar{u} }{|\nabla \bar{u}|}  \right) \ \ \hbox{in} \ \ \Omega \setminus \{x: |\nabla \bar{u}(x)|=0, \ \ a(x)\neq 0\}.
\end{equation}
\end{theorem}

Note that convergence of the Algorithm 1 does not require uniqueness of the minimizers of the problem (\ref{GMinProb}). In the case of the hybrid inverse problem where $a(x)$ is the magnitude of the current density vector field $J$ for some unknown conductivity $\sigma$ we can say more.

\begin{theorem}\label{theomainp}
Let $\Omega$ be a bounded region in $\subset \R^n$ with connected $C^{1,\alpha}$ bounday. Assume that $a=|J|>0$ a.e. in $\Omega$, where $J\in (L^2(\Omega))^n $ is the current density vector field  generated by an unknown conductivity $\sigma \in C^{\alpha}(\Omega)$ by imposing the voltage $f$ on $\partial \Omega$. Then the corresponding voltage potential $\bar{u}$ is the unique minimizer of $(\ref{GMinProb})$, and the sequences $b^k$, $d^k+\nabla u_f$, and $u^k+u_f$ produced by Algorithm 1 converge weakly to  $-J/ \lambda$, $\nabla \bar{u}$, and $\bar{u}$, respectively. 
\end{theorem}

\begin{remark} A generalization of the above theorem holds when $\Omega$ contains perfectly conducting $U_C$ or insulating $U_I$ inclusions \cite{MNT}. The sequences $b^k$, $d^k+\nabla u_f$, and $u^k+u_f$ produced by Algorithm 1 will converge weakly in $\Omega \setminus  U_C $ to  $-J/ \lambda$, $\nabla \bar{u}$, and $\bar{u}$, respectively.
\end{remark}

Our approach to (\ref{oldP}) is to start by working on the dual problem, which is well suited to a Douglas-Rachford splitting algorithm. This leads naturally to the alternating split Bregman algorithm. Indeed Setzer (\cite{Setzer, Setzer1}) proved the equivalence of the two methods.
In the process, we discovered that, for the inverse conductivity problem of interest here, the dual problem has a unique solution, which is in fact (up to a constant) the current $J$ ! (See \ref{niceCorl}).This shows that the alternating split Bregman algorithm is quite natural for the hybrid problem of Current Density Impedance Imaging.

The paper is organized as follows. In Section 2 we study the dual problem. Based on the analysis of the dual problem, in Section 3 we present our convergence results and prove Theorems \ref{theoMain} and \ref{theomainp}.  In Section 4 we will prove that the Algorithm 1 is stable with respect to possible errors in solving the Poisson equations (\ref{laplace}) at each step.  In section 5, we present several numerical experiments to demonstrate the convergence behaviour of our algorithm, in particular when the boundary function is not two-to-one.

\section{The dual problem}
Rockafellar-Fenchel duality \cite{Rb} will be the starting point for our proof of Theorem \ref{theoMain}. In this section we study the relation between the problem (\ref{GMinProb}) and its dual problem. To begin fix $u_f\in H^1(\Omega)$ with $u|_{\Omega}=f$ and let $a$ be a non-negative function in $L^2(\Omega)$. Define $F:(L^2(\Omega))^n\rightarrow \R$ and $G: H^1_0(\Omega)\rightarrow \R$ as follows
\begin{equation}\label{FandG}
F(d):=\int_{\Omega} a|d+\nabla u_f|dx, \ \ G(u)\equiv 0.
\end{equation}
Then the problem (\ref{GMinProb}) can be written as 
\[(P) \ \ \ \ \min_{u\in H^{1}_{0}(\Omega)} G(u)+F(\nabla u).\]
By Fenchel duality \cite{Rb}, the dual problem corresponding to the problem $(P)$ can be written as 
\begin{equation*}
(D) \hspace{.5cm}-\min_{b \in (L^2(\Omega))^n} \{ G^*(-\nabla \cdot b)+F^*(b)\},
\end{equation*}
Recall that the Legendre-Fenchel transform $F^*$ of a functional F on a Hilbert space $H$ is the function on $H^*$ given by
\[F^*(b)=\sup  \{\langle d,b \rangle -F(d): \ \  d\in H\}.\]

One can easily compute $G^*: H^{-1}(\Omega)\rightarrow \R$:
\begin{equation*}G^*(u)=\left\{ \begin{array}{ll}
0\ \ &\hbox{if} \ \ u\equiv 0,\\
\infty  \ \ &\hbox{if}  \ \ u \not \equiv 0. 
\end{array} \right.
\end{equation*} 
The following lemma provides a formula for $F^*:(L^2(\Omega))^n\rightarrow \R$. 

\begin{lemma} \label{f^*lem} Let $a\in L^2(\Omega)$ be non-negative and $u_f \in H^1(\Omega)$ with $u_f|_{\Omega}=f$. Then 

\begin{equation}\label{f^*}
F^*(b)=\left\{ \begin{array}{ll}
-\langle \nabla u_f, b \rangle\ \ &\hbox{if} \ \ |b(x)| \leq a(x)\ \ \hbox{a.e.} \ \ \hbox{in} \ \ \Omega,\\
\infty  \ \ & \hbox{otherwise}.
\end{array} \right.
\end{equation}
\end{lemma}
{\bf Proof.} Assume  
\[ |b(x)|> a(x),\]
on a subset $U$ of $\Omega$ with positive Lebesgue measure. Then 
\begin{eqnarray*}
F^*(b)&=&\sup_{d\in (L^2(\Omega))^n} \langle d, b \rangle-\int_{\Omega}a|d+\nabla u_f|\\
&=& -\langle b, \nabla u_f\rangle +\sup_{d\in (L^2(\Omega))^n} \left( \langle d, b \rangle-\int_{\Omega}a|d| dx\right)\\
&\geq &  -\langle b, \nabla u_f\rangle +\sup_{\lambda \in \R} \lambda \int_{U}(|b|^2-a(x)|b|)dx=\infty,\\
\end{eqnarray*}

where for the last inequality we choose $d(x)=\lambda b(x)$ for $x \in U$, and $d(x)=0$ otherwise.
Now assume 
\[|b(x)| \leq a(x), \ \ \hbox{a.e.}\]
then 
\begin{eqnarray}\label{f^**}
F^*(b)&=&-\langle b, \nabla u_f\rangle +\sup_{d\in (L^2(\Omega))^n} \left( \langle d, b \rangle-\int_{\Omega}a|d|dx \right) \\
&= &  -\langle b, \nabla u_f\rangle +\sup_{d\in (L^2(\Omega))^n} \int_{\Omega } (b\cdot d-a|d|)dx \nonumber  \\
&\leq &  -\langle b, \nabla u_f\rangle +\sup_{d\in (L^2(\Omega))^n} \int_{\Omega } |d(x)|(|b(x)|-a(x))dx \nonumber  \\
&\leq &  -\langle b, \nabla u_f\rangle .\nonumber 
\end{eqnarray}
Finally note that taking $d\equiv 0$ in (\ref{f^**}) yields 
\[F^*(b) \geq   -\langle b, \nabla u_f\rangle.\] \hfill $\Box$\\

It follows from (\ref{f^*}) that the dual problem can be explicitly written as
\begin{equation}\label{expDual} 
\max \{<\nabla u_f, b>: \ \ b\in (L^2(\Omega))^n,\ \ |b(x)|\leq a(x) \ \ a.e. \ \ \hbox{and} \ \ \nabla \cdot b\equiv 0\}.
\end{equation}
\\
Now let $u\in H^1_0(\Omega)$, $b \in (L^2(\Omega))^n$ with $\nabla \cdot b \equiv 0$, and $|b|\leq a$ a.e. in $\Omega$. Then
\begin{equation}\label{strongDual}
\int_{\Omega} a|\nabla u+\nabla u_f|dx\geq \int_{\Omega} |b||\nabla u+\nabla u_f|dx \geq \int_{\Omega}b \cdot (\nabla u+\nabla u_f)dx =<b, \nabla u_f>.
\end{equation}

Hence if we denote the optimal values of the primal and dual problem by $v(P)$ and $v(D)$, respectively, then $v(P)\geq v(D)$. This is a general fact. Moreover, since both of the functions $F$ and $G$ are continuous, it follows from Theorem 4.1 in \cite{ET} that strong duality holds, i.e.,  $v(P)=v(D)$ and the dual problem $(D)$ has an optimal solution.

The algorithm we propose seeks to construct a solution of the dual problem (D). In the rest of this section we show how this leads to the solution of the primal problem (P) as well. 

\begin{lemma}\label{pf^*}
Let $F^*$ be defined as (\ref{f^*}). Then
\small
{
\begin{equation*}\partial F^*(b)=\left\{ \begin{array}{ll}
-\nabla u_f+\{p: p \in (L^2(\Omega))^n, \ \ p=mb \ \ \hbox{on $\{x: a(x)>0\}$}, m \in M_b\} \ \ &\hbox{if} \ \  |b|\leq a \ \ a.e.\\
\emptyset \ \ & \hbox{otherwise},
\end{array} \right.
\end{equation*} 
}
where 
\[M_b:=\{m:\Omega \rightarrow [0,\infty ) \ \ |  \ \ \hbox{$m$ is measurable and } m(x)=0 \ \ \hbox{if} \ \ |b(x)|<a(x)  \}.\]
\end{lemma}
{\bf Proof.}  Assume $|b|\leq a$ a.e. in $\Omega$. Let  $m\in M_b$ and consider $d\in (L^2(\Omega))^n$ with $|d|\leq a$ a.e in $\Omega$. We have
\begin{eqnarray*}
\langle mb, d\rangle\leq \int_{\Omega} m|d||b|dx=\int_{\{ |b|=a\} } m|d||b|dx=\int_{\{ |b|=a \}}m|b|^2dx=\int_{\Omega}m|b|^2dx=\langle mb, b \rangle,
\end{eqnarray*}

and hence 
\[F^*(d)-F^*(b)=-\langle \nabla u_f, d \rangle + \langle \nabla u_f,b \rangle \geq  \langle -\nabla u_f+mb, d-b\rangle.\]
Note that the points where $a(x)=0$ (hence also $b(x)=d(x)=0$) do not contribute to any of the terms above. On the other hand 
\[F^*(d)-F^*(b) \geq  \langle -\nabla u_f+mb, d-b\rangle,\]
trivially holds if $|d|>a$ on a set of positive measure, as $F^*(d)=\infty$ in that case. Therefore   
\[-\nabla u_f+mb \in \partial F^*(b).\]
Now let $p\in \partial F^*(b)$ and define $\bar{p}:=\nabla u_f+p$.  Then for any $d\in (L^2(\Omega))^n$ with $|d|\leq a$ a.e in $\Omega$,
\[F^*(d)-F^*(b)=-\langle \nabla u_f, d \rangle + \langle \nabla u_f,b \rangle \geq  \langle p, d-b\rangle.\]
Consequently 
\begin{equation}\label{bd}
\langle \bar{p}, d\rangle \leq \langle \bar{p}, b\rangle,
\end{equation}
for all $d\in (L^2(\Omega))^n$ with $|d(x)|\leq a(x)$. In particular if we let 
\begin{equation*}d=\left\{ \begin{array}{ll}
a \frac{\bar{p}}{|\bar{p}|} \ \ & |\bar{p}|\neq 0\\
0 \ \ & \hbox{otherwise},
\end{array} \right.
\end{equation*} 
then it follows from (\ref{bd}) that 
 \begin{equation}\label{bd1}
  \langle \bar{p}, b\rangle \geq \langle\bar{p}, d \rangle =\int_{\Omega} a|\bar{p}| dx \geq \int_{\Omega} |b||\bar{p}| dx \ge \int_{\Omega} b\cdot \bar{p} \geq \langle \bar{p}, b\rangle.
 \end{equation}
Therefore all inequalities in (\ref{bd1}) are equalities. Hence there exists a non-negative function $m$ such that
 \begin{equation*}\bar{p}(x)=\left\{ \begin{array}{ll}
0\ \ &  \hbox{if}\ \ |b(x)|<a(x)\\
m b \ \ & \hbox{if} \ \ |b(x)|=a(x)\neq 0.
\end{array} \right.
\end{equation*} 
This completes the proof. \hfill $\Box$\\

Now we are ready to prove the following proposition which is a special case of the Rockafellar-Fenchel duality theorem \cite{Rb} in convex analysis. Given one solution of the dual problem (D), this result gives a description of all the solutions of the primal problem (P).

\begin{proposition} \label{ExpRF} Let $\hat{b}$ be an optimal solution of the dual problem (D). Then $\hat{u} \in H^1_0(\Omega)$ is an optimal solution of the primal problem (P) if and only if
\begin{equation}\label{heart}
\nabla \hat{u} \in \partial F^*(\hat{b}). 
\end{equation} 
\end{proposition}
{\bf Proof.} Let $\hat{u}$ be a solution of the primal problem.  Then, as we saw in (\ref{strongDual})
\begin{eqnarray*}
\int_{\Omega} a|\nabla \hat{u}+\nabla u_f|dx &\geq &  \int_{\Omega} |\hat{b}||\nabla \hat{u}+\nabla u_f|dx \\ 
&\geq &\int_{\Omega} \hat{b} \cdot (\nabla \hat{u}+\nabla u_f)dx = \langle \hat{b}, \nabla u_f\rangle.
\end{eqnarray*}
Since the duality gap is zero, both inequalities are equalities and
\begin{equation}\label{1b}
\hat{b}(x)=
a(x)\frac{\nabla \hat{u}+\nabla u_f}{|\nabla \hat{u}+\nabla u_f|} \ \ \hbox{if} \ \ |\nabla u+\nabla u_f|\neq 0.
\end{equation}
Therefore for $x$ with $a(x)\neq 0$, 
\[ \nabla \hat{u}(x)=-\nabla u_f+m(x)\hat{b}(x), \]
where 

 \begin{equation*}m(x)=\left\{ \begin{array}{ll}
\frac{ |\nabla u(x)+\nabla u_f(x)|}{a(x)}\ \ &  \hbox{if}\ \  |\nabla u+\nabla u_f|\neq 0\\
0 \ & \hbox{otherwise}.
\end{array} \right.
\end{equation*} 
In view of Lemma \ref{pf^*}, we conclude $\nabla \hat{u} \in \partial F^*(\hat{b})$.  

Now assume $\nabla \hat{u} \in F^*(\hat{b})$. Then by Lemma \ref{pf^*}, there exists $m \in M_{\hat{b}}$ such that for $x$ with $a(x)\neq 0$,
$ \nabla \hat{u}(x)=-\nabla u_f+m(x)\hat{b}(x)$.  Therefore
\[ \langle \hat{b}, \nabla u_f\rangle =\int_{\Omega} \hat{b} \cdot (\nabla \hat{u}+\nabla u_f)dx =\int_{\Omega} |\hat{b}||\nabla \hat{u}+\nabla u_f|dx=\int_{\Omega} a|\nabla \hat{u}+\nabla u_f|dx,\]
which means $\hat{u}$ is a minimizer of the primal problem (P). \hfill $\Box$\\

In the next section, we shall see that the above relation (\ref{heart}) between optimal solutions of the primal and dual problems is at the heart of Algorithm 1.  The relation (\ref{1b}) shows that $\hat{b}$ is determined on the set where $|\nabla \hat{u}(x)+\nabla u_f(x)|$ does not vanish. We record this fact as a separate partial uniqueness result:
\begin{proposition}\label{uniqDualProSol}
Let $\hat{u}$ be an optimal solution of the primal problem and assume $\hat{b}_1$ and $\hat{b}_1$ are two optimal solutions for the dual problem (D). Then 
\begin{equation}\label{b_1=b_2}
\hat{b}_1\equiv \hat{b}_2 \ \ \hbox{in} \ \ \Omega \setminus \{x:  |\nabla \hat{u}(x)+\nabla u_f(x)|=0, a(x)\neq 0\}.
\end{equation}

\end{proposition}
{\bf Proof.} By Lemma \ref{ExpRF} $\nabla \hat{u} \in \partial F^*(\hat{b}_1) \cap \partial F^*(\hat{b}_1)$. hence it follows from lemma \ref{pf^*} that there exist $m_1 \in M_{\hat{b}_1}$ and $m_2 \in M_{\hat{b}_1}$ such that 
\[\nabla \hat{u}=-\nabla u_f+m_1\hat{b}_1 \ \ \hbox{and} \ \ \nabla \hat{u}=-\nabla u_f+m_2\hat{b}_2. \]
Thus $m_1\hat{b}_1\equiv m_1\hat{b}_2$. Since $\hat{b}_1$ and $\hat{b}_1$ are both optimal, by (\ref{1b})
\[|\hat{b}_1(x)|=|\hat{b}_2(x)|=a(x) \ \ \hbox{on} \ \ \Omega \setminus \{x:  |\nabla \hat{u}(x)+\nabla u_f(x)|=0, \ \ a(x)\neq 0\}.\]
Hence (\ref{b_1=b_2}) follows. \hfill $\Box$

If $a(x)$ is the magnitude of the current corresponding to a conductivity $\sigma$ the above proposition yields uniqueness of solutions to the dual problem. 

\begin{corollary}\label{niceCorl}
Let $\Omega$ be a bounded region in $\subset \R^n$ with connected $C^{1,\alpha}$ bounday. Assume $a=|J|>0$ a.e. in $\Omega$, where $J \in (C^{\alpha}(\overline{\Omega}))^n$ is the current density vector field  generated by an unknown conductivity $\sigma \in C^{\alpha}(\Omega)$ by imposing the voltage $f$ on $\partial \Omega$. Then the corresponding voltage potential $\bar{u}$ is the unique minimizer of $(\ref{GMinProb})$ and the dual problem $(D)$ has a unique solution $\hat{b}$. Furthermore $\hat{b}= -J$. 
\end{corollary}
{\bf Proof.} The proof follows from Proposition \ref{uniqDualProSol} , the uniqueness Theorem 1.3 in (\cite{NTT08}) and equation (\ref{1b}). \hfill $\Box$

\begin{remark} The above corollary generalizes to the case when $\Omega$ contains perfectly conducting ($U_P$) and/or perfectly insulating $U_I$ inclusions. Proposition \ref{uniqDualProSol} together with the uniqueness result in \cite{MNT} show that a solution $\hat{b}$ of the dual problem (D) equals the current $J$ in $\Omega \setminus U_p$. We omit the details.
\end{remark}

\section{Convergence analysis}
In this section we present proofs of Theorems \ref{theoMain} and \ref{theomainp}. The proofs rely on the representation of solutions of the primal problem (\ref{heart}) in Proposition \ref{ExpRF}.  Notice that the dual problem (D) can be written in the form of an inclusion problem
 \begin{equation}\label{inclusion}
 0\in A(\hat{b})+B(\hat{b}),
 \end{equation}
 where $A:=\partial (G^*o(-\nabla ^*))$ and $B:=\partial F^*$ are maximal monotone operators on $(L^2(\Omega))^n$. If we can compute a solution $\hat{b}$  of the problem (\ref{inclusion}) as well as $\hat{d}\in B(\hat{b})=\partial F^*(\hat{b})$ then, by Proposition \ref{ExpRF},  $\hat{u}= \nabla^{-1}(\hat{d})$ will be a solution of the primal problem (P).  
The Douglas-Rachford splitting  method in Convex Analysis yields precisely such a pair $(\hat{b},\hat{d})$. Following this route to the primal problem leads naturally to the Alternating Split Bregman algorithm, as we explain below.

Let $H$ be a real Hilbert space and let $A, B: H\rightarrow 2^{H}$ be two maximal monotone (set valued) operators. For a set valued function $P: H\rightarrow 2^{H}$, let $J_P$ denote its resolvent i.e.,
\[J_{P}=(Id+P)^{-1}.\]
The sub-gradient of a convex, proper, lower semi-continuous function is maximal monotone and if $P$ is maximal monotone then $J_P$ is single valued \cite{FPA, Rb} .  Lions and Mercier \cite{LM}  showed that for any general maximal monotone operators $A,B$ and any initial element $x_0$ the sequence defined by the Douglas-Rachford recursion:
\begin{equation}\label{iter:DR}
x_{k+1}=(J_A(2J_B-Id)+Id-J_B)x_k,
\end{equation}
converges weakly to some point $\hat{x} \in \emph{H}$ such that $\hat{p}=J_B (\hat{x})$ solves the inclusion problem (\ref{inclusion}). Recently Svaiter \cite{S} proved that the sequence $p_k=J_{\eta B}(x_k)$ also converges weakly to $\hat{p}$. This fact will be important for our problem. The following theorem describes the Douglas-Rachford splitting  algorithm and summarizes these known convergence results.

\begin{theorem}\label{DRS}
Let $H$ be a Hilbert space and let $A, B: \emph{H}\rightarrow 2^{\emph{H}}$ be maximal monotone operators and assume that a solution of (\ref{inclusion}) exists. Then, for any initial elements $x_0$ and $p_0$ and any $\lambda>0$, the sequences $p_{k}$ and $x_k$ generated by the following algorithm
\begin{eqnarray}
x_{k+1}&=&J_{\lambda A}(2p_k-x_k)+x_k-p_k \nonumber \\
p_{k+1}&=&J_{\lambda B}(x_{k+1}),
\end{eqnarray}
converge weakly to some $\hat{x}$ and $\hat{p}$ respectively. Furthermore, $\hat{p}=J_{\lambda B}(\hat{x})$  and $\hat{p}$ satisfies
\[0\in A(\hat{p})+B(\hat{p}).\]
\end{theorem}

We wish to apply the Douglas-Rachford splitting algorithm to the operators $A:=\partial (G^*o(-\nabla ^*))$ and $B:=\partial F^*$. We need an efficient way to evaluate the resolvents $J_{\lambda A}(2p_k-x_k)$ and $J_{\lambda B}(x_{k+1})$ at each iteration. If we let 
\begin{equation}
x_k=\lambda (b^k+d^k), \ \ \ \ p_k=\lambda  b^k, \ \ k\geq 0, \\
\end{equation}
then to evaluate $J_{\lambda A}(2p_k-x_k)$ and $J_{\lambda B}(x_{k+1})$  for all $k\geq 0$,  it is enough to find the minimizers $u^{k+1}$ and $d^{k+1}$ of the functionals 
\begin{equation}\label{Iu1}
I^k_1(u)=\parallel \nabla u+b^k-d^k\parallel ^2,
\end{equation}
and
\begin{equation}\label{Iu2}
I^k_2(d)=\int_{\Omega} a|d+\nabla u_f|dx+\parallel b^k+\nabla u^{k+1}-d\parallel^2
\end{equation}
and set $b^{k+1}=b^k+\nabla u^{k+1}-d^{k+1}$. Indeed the resolvents $J_{\lambda A}(2p_k-x_k)$ and $ J_{\lambda B}(x_{k+1})$ can be computed as follows 
\[J_{\lambda A}(2p_k-x_k)=\lambda(b^k+\nabla u^{k+1}-d^k),\]
and 
\[J_{\lambda B}(x_{k+1})=\lambda (b^k+\nabla u^{k+1}-d^{k+1}),\]
 (see \cite{Setzer, Setzer1} for a proof).  Finding the minimizer of (\ref{Iu1}) amounts to solving a Poisson equation 
\[\Delta u^{k+1}=\nabla \cdot (d^k-b^k), \ \ u^{k+1}|_{\partial \Omega}=0.\] 
 As well, the minimizer of the functional $I^k_2(d)$ can be computed explicitly as follows

\[d^{k+1}:=\left\{ \begin{array}{ll}
\max \{|\nabla u^{k+1}+\nabla u_f +b^k|-\frac{a}{\lambda},0\}\frac{\nabla u^{k+1}+\nabla u_f +b^k}{|\nabla u^{k+1}+\nabla u_f  +b^k|} -\nabla u_f \ \ &\hbox{if} \ \ |\nabla u^{k+1}(x)+ \nabla u_f +b^k(x)|\neq 0,\\
-\nabla u_f   \ \ &\hbox{if} \ \ |\nabla u^{k+1}(x)+ \nabla u_f +b^k(x)|=0.
\end{array} \right.
\]

We are thus led to the Algorithm 1 for simultaneously finding solutions of both the primal problem and the dual problem. To prove the strong convergence of the series (\ref{conv-est}) we will need the following simple lemma on firmly non-expansive operators.

It is well known that the operator $T:=J_A(2J_B-Id)+Id-J_B$ is firmly non-expansive, i.e.,  $T=\frac{1}{2}Id+\frac{1}{2}R$ such that $R$ is non-expansive:
\[\parallel Rx -Ry \parallel \leq \parallel x-y\parallel \ \ \ \ \hbox{for all} \ \ x,y \in H.\]

\begin{lemma}\label{nonexpansice-lemma}  If $T:H\rightarrow H$ is a firmly non-expansive operator and $x_{k+1}=T(x_k)$ with $x_0 \in H$, then
\[ \parallel x_{k+1}-\hat{x}\parallel ^2+\parallel x_{k+1}-x_k\parallel ^2 \leq \parallel x_{k}-\hat {x} \parallel ^2.\]
\end{lemma}
{\bf Proof. } Since $T$ is firmly non-expansive, $R=2T-Id$ is a non-expansive operator.

 Hence
\[\parallel R x_k- R \hat{x}\parallel ^2 =-\parallel x_k- \hat{x} \parallel ^2+2\parallel T x_k -T\hat{x}\parallel^2-2 \parallel (Id-T)x_k -(Id-T)\hat{x}\parallel^2. \]
Therefore we have
\[\frac{1}{2}(\parallel x_k- \hat{x} \parallel ^2 -\parallel R x_k- R \hat{x}\parallel ^2 )=\parallel x_k- \hat{x} \parallel ^2-\parallel x_{k+1} -\hat{x}\parallel^2- \parallel x_{k+1}-x_k\parallel^2.\]
Since $R$ is non-expansive, the left hand side of the above inequality is non-negative, which completes the proof. \hfill $\Box$\\

{\bf Proof of Theorem \ref{theoMain}.} By interpreting Algorithm 1 as a Douglas-Rachford splitting algorithm as detailed above, weak convergence of the sequences  $d^k$, and $b^k$  follows immediately from Theorem \ref{DRS}. To prove the estimate (\ref{conv-est}), let $T=J_{\lambda A}(2J_{\lambda B}-Id)+Id-J_{\lambda B}$. Since $T$ is firmly non-expansive, by Lemma \ref{nonexpansice-lemma} we have
\begin{eqnarray}
\parallel x_{k+1}-\hat{x}\parallel ^2+\parallel x_{k+1}-x_n\parallel ^2 \leq \parallel x_{k}-\hat {x} \parallel ^2,
\end{eqnarray}
where $\hat{x}$ is the weak limit of $x_k$ with $T(\hat{x})=\hat{x}$.  By the above inequality we have
\begin{equation}
\sum_{k=0}^{\infty}\parallel x_{k+1}-x_{k}\parallel^2 <\infty. 
\end{equation}
Now observe that
\[x_{k+1}-x_k=\lambda ( b^{k+1}+d^{k+1}-b^k-d^{k})=\lambda (\nabla u^{k+1}-d^k),\]
and hence (\ref{conv-est}) follows.  

By Theorem \ref{DRS},  $\hat{p}=\lambda \hat{b}$ is a minimizer of the dual problem and $J_{\lambda \partial F^*}  ( \lambda (\hat{d}+\hat{b}))=\lambda \hat{b}$.  Therefore
\[\lambda \hat{b}+\lambda \partial F^* (\lambda \hat{b})=\lambda (\hat{d}+\hat{b}) \ \ \Leftrightarrow \ \ \hat{d} \in \partial F^*(\lambda \hat{b}).\]

In view of (\ref{conv-est}) the sequence $\{u^k\}_{k\in N}$ is bounded in $H^1_0(\Omega)$ and therefore it has a weakly converging subsequence. Let $\hat{u}$ be a weak cluster point of the sequence $\{u^k\}$.  Then $\nabla\hat{u} = \hat{d} \in \partial f^*(\hat{p})$ and hence by Proposition \ref{ExpRF} $\hat{u}$ is a solution of the primal problem (P). On the other hand $\nabla\hat{u}=\hat{d}$ for every weak cluster point $\hat{u}$ of the sequence $\{u^k\}_{k\in N}$. Since $\nabla$ is injective on $H^1_0(\Omega)$,  $\{u^k\}_{k\in N}$ has at most one weak cluster point. In view of Proposition \ref{uniqDualProSol}, the proof is now complete. \hfill $\Box$\\

{\bf Proof of Theorem \ref{theomainp}.} The proof follows from \ref{theoMain}, combined with the uniqueness Theorem 1.3 in \cite{NTT08}, and our Corollary \ref{niceCorl}. \hfill $\Box$

\section{Approximate alternating split Bregman algorithm }

In this section we show that the alternating split Bregman algorithm converges to the correct solutions even in the presence of possible errors at each step in solving Poisson  equations.  The proof relies on the following theorem about the Douglas-Rachford splitting algorithm.  

\begin{theorem} \label{theo:SP}{\bf (Svaiter \cite{S})} Let $\lambda>0$, and let $\{\alpha_k\}_{k \in N}$ and $\{\beta_k\}_{k \in N}$ be sequences in a Hilbert space $H$. Suppose $0 \in\hbox{ran}(A+B)$, and $\sum_{k \in N} (\parallel \alpha _k\parallel +\parallel \beta_k\parallel)<\infty$.  Take $x_0 \in H$ and set
\begin{equation}
x_{k+1}=x_k+J_{\gamma A}(2(J_{\lambda B}x_k+\beta_n)-x_k)+\alpha_k -(J_{\lambda B}x_k+\beta_k) , \ \ k\geq 1. 
\end{equation}
Then $x_k$ and $p_k=J_{\lambda B}x_k$ converge weakly to $\hat{x} \in H$ and  $\hat{p} \in H$, respectively and $\hat{p}=J_{\lambda B}\hat{x} \in (A+B)^{-1}(0)$. \\
\end{theorem}

The proof of the above theorem in infinite-dimensional Hilbert spaces is due to Svaiter \cite{S} (see also \cite{C}). It suggests the following approximate version of our algorithm\\

\textbf{Approximate alternating split Bregman algorithm:} \\

Let $u_f \in H^1(\Omega)$ with $u_f|_{\partial \Omega}=f$  and initialize $b^0, d^0 \in (L^2(\Omega))^n$. For $k\geq 1$\\ 
\begin{enumerate}
\item Find an approximate solution $u^k$  of 
\[ \Delta u^{k+1}= \nabla \cdot (d^k(x)-b^k(x)), \ \ u^{k+1}|_{\partial
\Omega}=0,\]
with $\parallel \nabla u^k-\nabla u^k_{ex}\parallel \leq \alpha_k$, where $u^{k+1}_{ex}$ is the exact solution of the above problem. 
\item Compute
 \small{
\[d^{k+1}:=\left\{ \begin{array}{ll}
\max \{|\nabla u^{k+1}+\nabla u_f +b^k|-\frac{a}{\lambda},0\}\frac{\nabla u^{k+1}+\nabla u_f +b^k}{|\nabla u^{k+1}+\nabla u_f  +b^k|} -\nabla u_f  &\hbox{if} \ \ |\nabla u^{k+1}(x) +\nabla u_f +b^k(x)|\neq 0,\\
-\nabla u_f    &\hbox{if} \ \ |\nabla u^{k+1}(x) +\nabla u_f+ b^k(x)|=0.
\end{array} \right.
\] 
 }
\item Set
\[b^{k+1}(x):=b^{k}(x)+\nabla u^{k+1}(x)-d^{k+1}(x).\]
\end{enumerate}

 By Theorem \ref{theo:SP} and an argument similar to that of Theorem \ref{theoMain}  we can prove the following theorem about convergence of the sequences $u^k$, $d^k$, and $b^k$ produced by the above algorithm

\begin{theorem}\label{theo2}
Let $a \in L^2(\Omega)$ be a non-negative function and $f\in H^{1/2}(\Omega)$. If

\[\sum_{k=1}^{\infty} \alpha_k<\infty,\]

then for any $u_f \in H^1(\Omega)$ and any $b^0, d^0 \in (L^2(\Omega))^n$ the sequences $\{b^k\}_{k\in N}$, $\{d^k\}_{k\in N}$, and $\{u^k\}_{k\in N}$ produced by the approximate alternating split Bregman algorithm converge weakly to some $\hat{b}$, $\hat{d}$, and $\hat{u}$. Moreover  $\{\nabla u^{k+1}-d^k\}_{k\in N}$ converges strongly to zero and
\begin{equation*}
\sum_{k=0}^{\infty}\parallel\nabla u^{k+1}-d^k \parallel ^2<\infty.
\end{equation*}
Furthermore $\bar{u}:=\hat{u}+u_f$ is a solution of the minimization problem ($\ref{GMinProb}$), $\hat{b}$ is a solution of the dual problem (D),  $\hat{d}=\nabla \hat{u}$, $\nabla \cdot \hat{b} \equiv 0$, and
\begin{equation}\label{bhat}
\hat{b}=\frac{1}{\lambda} \left( a\frac{\nabla \bar{u} }{|\nabla \bar{u}|}  \right) \ \ \hbox{in} \ \ \Omega \setminus \{x: |\nabla \bar{u}(x)|=0, \ \ a(x)\neq 0\}.
\end{equation}
\end{theorem}

There is also an analogue to Theorem \ref{theomainp} for the Approximate alternating split Bregman algorithm, which we omit.

\section{Numerical experiments}
In this section we study numerically the convergence behavior of the proposed alternating split Bregman algorithm.  In a model problem, we cseek to reconstruct the conductivity from knowledge of the magnitude of one current density $|J|$ given inside the unit square $\Omega = (0,1)\times (0,1)$ and  the corresponding voltage potential $f$ on  $\partial \Omega$. In \cite{NTT08} authors presented a simple iterative algorithm to recover the conductivity $\sigma$ from the knowledge of $(|J|,f)$. To be well defined, this algorithm requires the solution $\hat{u}$ and all intermediate functions $u^k$ produced by the algorithm to have non-vanishing gradients in $\Omega$. The alternating split Bregman algorithm proposed in this paper does not require $|\nabla \hat{u}|>0$ in $\Omega$ and converges (by Theorem \ref{theoMain}) to an optimal solution of (\ref{GMinProb}). We perform several numerical experiments to illustrate this convergence behavior. 

\subsection{Data simulation}
\begin{figure}[!h]
  \centering
  \includegraphics[height=70mm]{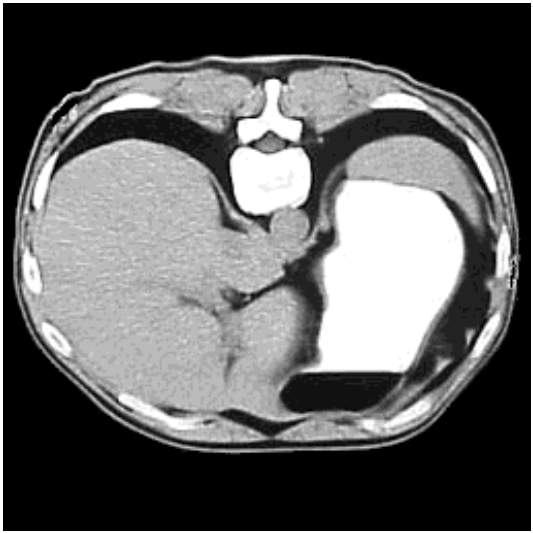}
  \caption{The original conductivity distribution used in the data simulation. }
  \label{figure1}
\end{figure}

To simulate the internal data $|J|$, we use an abdominal human CT image rescaled to a realistic range of tissue conductivity, with values varying from 1 to 1.8 S/m The scaled conductivity distribution, on a uniform grid $128\times 128$, is shown in Figure \ref{figure1}.

Given the above conductivity distribution, we first solve numerically the Dirichlet problem
\[
\nabla\cdot\sigma\nabla v = 0, \ \ v|{\partial\Omega} = f,
\]
on the grid $128\times 128$. To provide high accuracy of the numerical solution, we look for a solution of the form $v = u_{h} + u$, where $u_h$ is the harmonic function satisfying the Dirichlet boundary condition and $u$ is the solution to the Poisson equation
\[
- \nabla\cdot\sigma\nabla u = \nabla\cdot\sigma\nabla u_h
\]
with the zero boundary conditions. For the range of conductivity values described above, the norm $||u||$ is small compared to that of $||u_h||$. Such a representation is also helpful in the Algorithm 1; the function $u_h$ needs to be computed just once. As forward solvers, we use the FORTRAN sofware FISHPACK and MUDPACK for elliptic problems (see \cite{fortran}). We run those routines with the double precision. Comparison of the numerical solutions with an analytical one (in the case of constant conductivity) verified that the relative $L^{2}$-error does not exceed $10^{-6}$. Note that if such an error were to exceed $10^{-5}$, then it may severely affect the quality of the reconstructed images. We combined the above solvers with the numerical differentiation via the three- or five-point Lagrangian interpolation, to preserve the high accuracy needed  in the reconstruction algorithms. Once the solution $u$ is computed the magnitude of the current density in $\Omega$, which is the data for our problem, is simulated as $|J| = \sigma |\nabla u|$.

\subsection{Reconstruction algorithms}

For reconstruction we use the alternating split Bregman algorithm, Algorithm 1, choosing for $u_f$ the harmonic extension $u_h$ of the boundary voltage $f$. This leads to a simpler algorithm. For comparison purpose, we will also apply the so-called simple iterations algorithm from \cite{NTT08}. Below we outline both of algorithms for convenience. \\

{\bf Algorithm 2 (Simplified Alternating split Bregman algorithm)}\\ \\
Let $u_h$ be the harmonic extension of $f$ to $\Omega$ and initialize $d^0,b^0 \in (L^2(\Omega))^n$. \\
\begin{enumerate} 
\item Beginning with $k=0$, solve 
\begin{equation}\label{laplace}
 \Delta u= \nabla \cdot (d^k(x)-b^k(x)), \ \ u|_{\partial
\Omega}=0, 
\end{equation}
and let $v^{k+1}=u+u_h$.
\item Compute
\small{
\begin{equation*}\label{expEq}
d^{k+1}(x):=\left\{ \begin{array}{ll}
\max \{|\nabla v^{k+1}(x) +b^k(x)|-\frac{a(x)}{\lambda},0\}\frac{\nabla v^{k+1}(x)+b^k(x)}{|\nabla v^{k+1}(x) +b^k(x)|} \ \ &\hbox{if} \ \ |\nabla v^{k+1}(x) +v^k(x)|\neq 0,\\
0  \ \ &\hbox{if} \ \ |\nabla v^{k+1}(x) +b^k(x)|=0.
\end{array} \right.
\end{equation*}
}
\item Let \[b^{k+1}(x)=b^{k}(x)+\nabla v^{k+1}(x)-d^{k+1}(x).\]
\item If $\frac{||v_{k+1} - u_{k}||}{||v_{k+1}||}\leq Tol$, then compute
\[
\sigma = \frac{|J|}{|\nabla{v}_{k+1}|}.
\]
Otherwise, $k := k + 1$ and repeat the process.\\
\end{enumerate}

{\bf The simple iterations algorithm.} \\

Let $u_h$ be the harmonic extension of $f$ in $\Omega$ and initialize $u_{0} = u_f, \sigma_{1} = \frac{|J|}{\nabla{u}_{0}}$. \\
\begin{enumerate}
\item Beginning with $k = 1$, solve the problems
\[
-\Delta u = \nabla\cdot(\sigma_{k}\nabla{u}_{h}),~~u_{\mid{\partial\Omega}} = 0. 
\]
$v_{k} = u + u_{h}$. \\
\item Update 
\[
\sigma_{k+1} = \frac{|J|}{|\nabla{v}_{k}|}.
\]
\item If $\frac{||\sigma_{k+1} - \sigma_{k}||}{||\sigma_{k+1}||}\leq Tol$, then STOP. Otherwise, let $k := k + 1$ and repeat the process.
\end{enumerate}

In both algorithms, the harmonic part $u_h$ of the solution is computed via the over-relaxation method, and the Dirichlet problem for the Poisson equation is solved numerically by using the implicit conjugate-gradient method (see \cite{samar}) in which the iterative process is constructed by minimizing the error of each approximation in the energy norm.  The correction vector from the Krylov space is determined on each iteration. By virtue of the implicit method, a five-diagonal matrix is inverted on each iteration using a preconditioner.        

\subsection{Numerical reconstructions}
\begin{figure}[h^*]
  \centering
  \includegraphics[height=80mm]{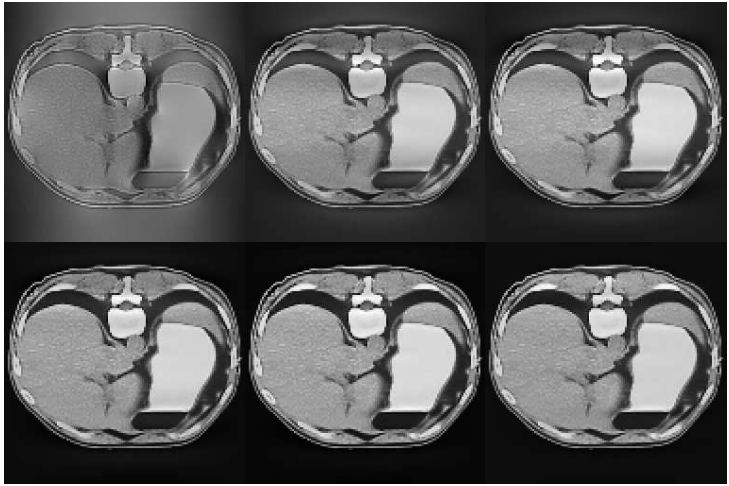}
  \caption{Conductivity reconstruction using the alternating split Bregman algorithm with $N = 1, 5, 10, 30, 50, 100$ iterations (shown from the left upper corner to the right lower corner) for the almost two-to-one boundary condition $f(x,y) = y$.}
  \label{figure3}
\end{figure}

\begin{figure}[!h]
  \centering
  \includegraphics[height=80mm]{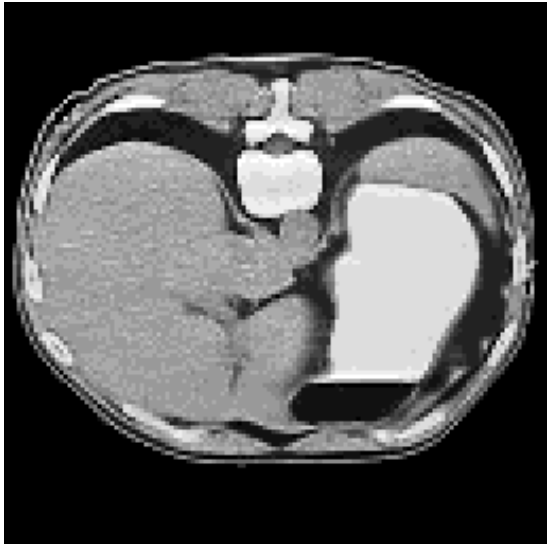}
  \caption{Conductivity reconstructed using the alternating split Bregman algorithm with $N=60$ iterations for the almost two-to-one boundary condition $f(x,y) = y$. }
  \label{figure3prime}
\end{figure}

We use $\lambda = 1$ in all the numerical experiments with the alternating Bregman algorithm. In our first experiment we choose the almost two-to-one boundary voltage $f(x,y) = y$. The results obtained by applying the split Bregman algorithm with $N = 1, 5, 10, 30, 50, 100$ iterations are shown in Figure \ref{figure3}. A larger image of  the conductivity reconstructed using the alternating split Bregman algorithm with $N=60$ is shown in Figure \ref{figure3prime}. This image may be compared with the original image in Figure \ref{figure1}.

\newpage
\begin{figure}[!h]
  \centering
  \includegraphics[height=80mm]{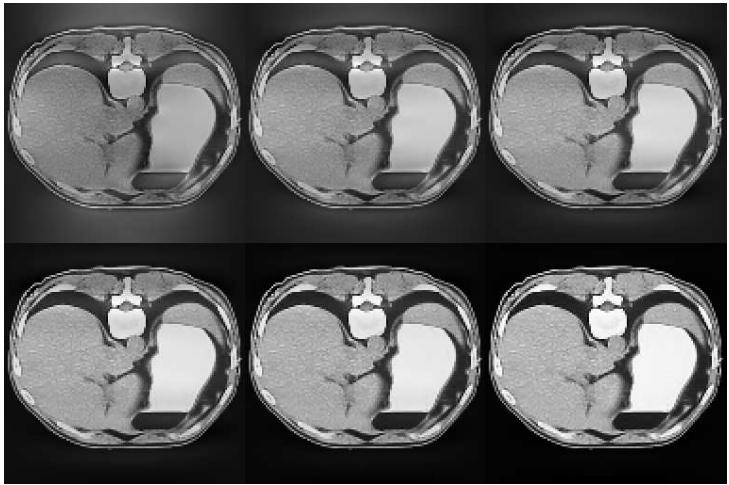}
  \caption{Reconstruction using the simple iterations algorithm with $N = 1, 5, 10, 30, 50, 100$ iterations (show from the left upper corner to the right lower corner) for the almost two-to-one boundary condition $f = y$.}
  \label{figure4}
\end{figure}

For comparison we repeat the above experiment with the same almost two-to-one boundary condition using the simple iterations algorithm. Figure \ref{figure4} shows the resulting conductivity for different number of iterations.  For two-to-one boundary data, the results of the two algorithms are similar, although the simple iterations method slightly outperforms the alternating Bregman algorithm in this case. \\ \\
The Tables \ref{table1} and \ref{table2} give the numerical errors of experiments with alternating split Bregman and simple iterations algorithms for different levels of tolerance and the computation times on a Dell Precision T5400 workstation with an Intel Xeon 64-bit 2 core processor.

\begin{table}[!h]
    \begin{tabular}{ | l | l | l | l |}
    \hline
    Tolerance & Relative $L^{2}-$ error vs EXACT &  Elapsed time(s)  & Total number of iterates    \\ \hline
    $5x10^{-5}$  &  0.0156  & 703  & 122 \\ \hline
     $1x10^{-4}$   & 0.0148 & 574  & 99 \\ \hline
    $2x10^{-4}$  & 0.0075 & 443 & 76  \\ \hline
   $5x10^{-4}$    &  0.0166   &   276 & 47  \\ \hline   
    \end{tabular}
   \caption{Numerical errors and elapsed times for alternating split Bregman}    
     \label{table1}    
   \end{table}

\begin{table}[!h]
    \begin{tabular}{ | l | l | l | l |}
    \hline
    Tolerance & Relative $L^{2}-$ error vs EXACT &  Elapsed time(s)  & Total number of iterates    \\ \hline
    $5x10^{-5}$  &  0.0030  & 672  & 110 \\ \hline
    $ 1x10^{-4}$   & 0.0030 & 583  & 99 \\ \hline
    $2x10^{-4}$  & 0.0137 & 446 & 73  \\ \hline
   $5x10^{-4}  $  & 0.0141   &   264 & 43  \\ \hline   
      \end{tabular}
    \caption{Numerical errors and elapsed times for simple iterations}
  \label{table2}                                     
\end{table}

\begin{figure}[!h]
  \centering
  \includegraphics[height=70mm]{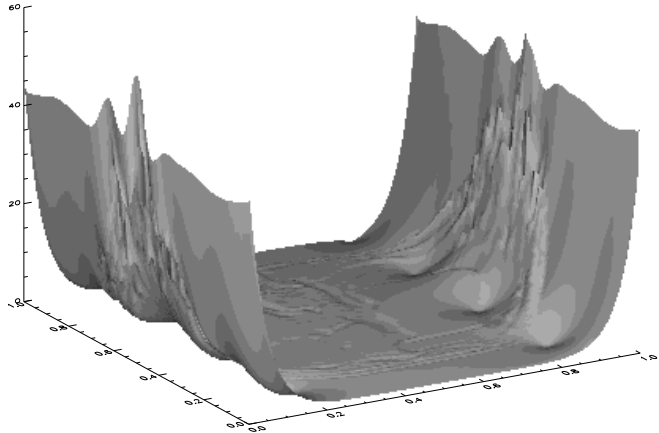}
  \caption{Magnitude of the current density $|J|$ for the non two-to-one boundary data $f(x,y)=y+2\sin(7 \pi y)$.}
  \label{figure4prime}
  \label{|J|}
\end{figure}

\begin{figure}[h*]
  \centering
  \includegraphics[height=80mm]{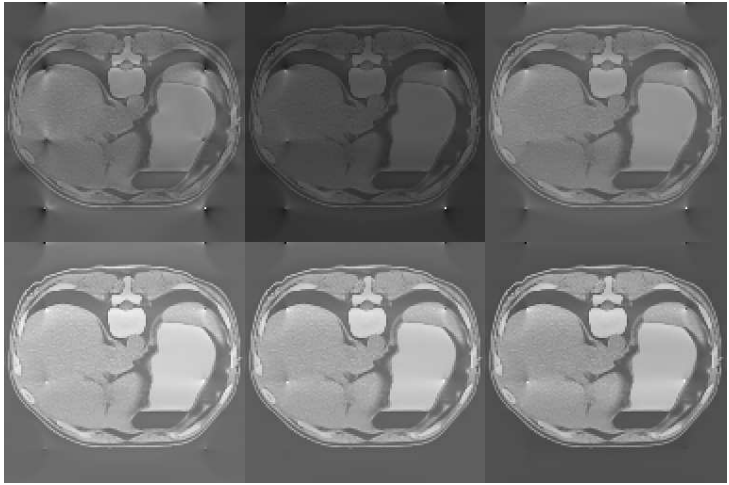}
  \caption{Conductivities constructed using the alternating split Bregman algorithm  with $N = 1, 5, 10, 30, 50, 100$ iterates (shown from the left upper corner to the right lower corner) for the non two-to-one boundary data $f(x,y)=y+2\sin(7 \pi y)$.}
  \label{figure5}
\end{figure}

Recall that the simple iterations algorithm requires $|\nabla u^k|>0$ on $\Omega$ for all $k \geq1$. In general if the boundary data $f$ is not two-to-one then there exist $x\in \Omega$ such that $\nabla \hat{u}(x)=0$. This may lead to divergence of the simple iterations algorithm. For the next experiments, we chose the boundary voltage $f(x,y)=y+2\sin(7 \pi y)$, which is not two-to-one. As shown in Figure \ref{|J|} the resulting surface $z=|J|$ touches the $xy$-plane. For this boundary data, the simple iterations algorithm breaks down. However as shown in Figure \ref{figure5} the alternating split Bregman algorithm converges.

\newpage

\begin{figure}[!h]
  \centering
  \includegraphics[height=70mm]{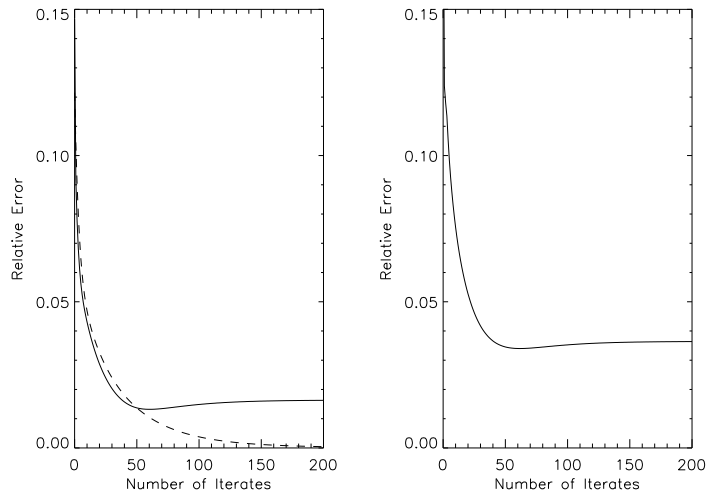}
  \caption{Rate of convergence for alternating split Bregman (Solid) and simple iterations (Dashed) for two-to-one boundary data $f(x,y)=y$ (left) and non two-to-one data $f(x,y)=y+2\sin(7 \pi y)$. }
  \label{rate}
\end{figure}

In Figure \ref{rate}, we show the rate of convergence of the alternating split Bregman (Solid) and simple iterations (Dashed) algorithms for two-to-one boundary data $f(x,y)=y$ (left) and non two-to-one data $f(x,y)=y+2\sin(7 \pi y)$. There is no dashed curve on the left, as the errors rapidly exceed the scale in the figure. Our numerical experiments indicate that simple iterations slightly outperform the alternating split Bregman algorithm for two-to-one boundary data. However, for non two-to-one boundary data simple iterations may diverge, while the alternating split Bregman always converges (see Figure \ref{rate}).

\begin{figure}[!h]
  \centering
  \includegraphics[height=80mm]{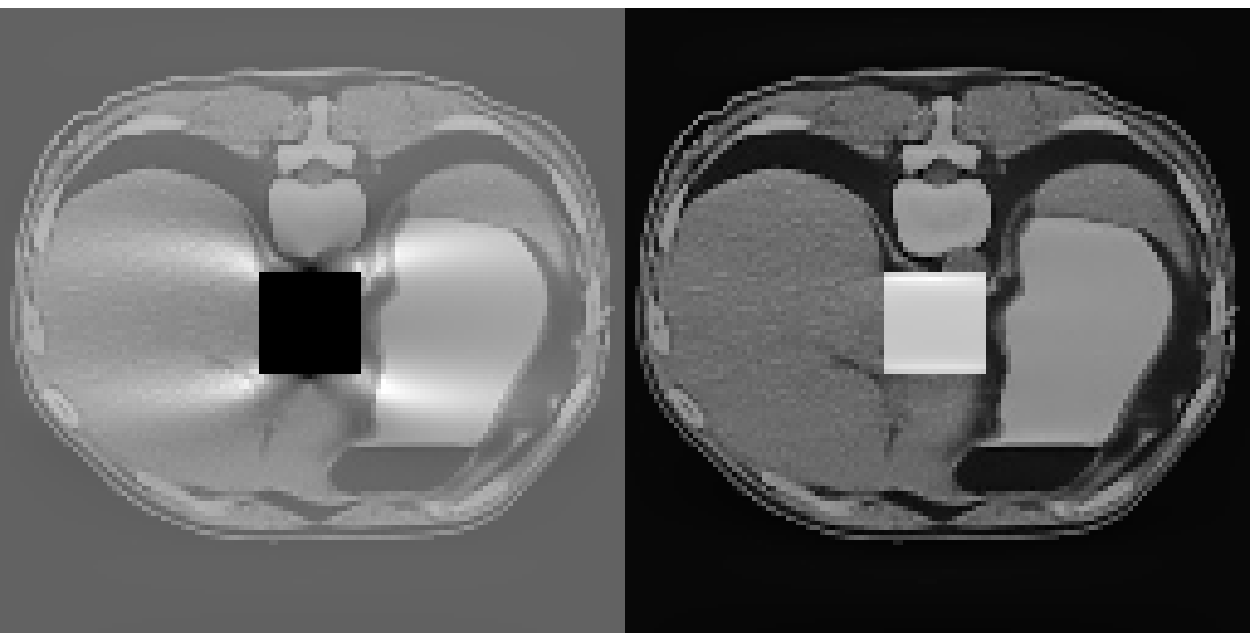}
  \caption{Reconstruction in the presence of the perfectly conducting (right) and insulating (left) inclusions  }
  \label{inclu}
\end{figure}

In \cite{MNT} the more genral problem of recovering an isotropic conductivity outside some perfectly conducting or insulating inclusions was considered. The data was, as before, the magnitude of one current density field $|J|$ in the interior of $\Omega$.. We proved that (except in some exceptional cases) the conductivity outside the inclusions, and the shape and position of the perfectly conducting and insulating inclusions are uniquely determined by the magnitude of the current generated by imposing a given boundary voltage. Since the relevant minimization problem is still of the form \ref{min_prob} in this case, the split Bregman algorithm can be applied. Figure \ref{inclu} shows the conductivity constructed using $100$ iterations of the alternating split Bregman algorithm in the presence of perfectly conducting (right) and insolating (left) inclusions. 

\begin{figure}[!h]
  \centering
  \includegraphics[height=100mm]{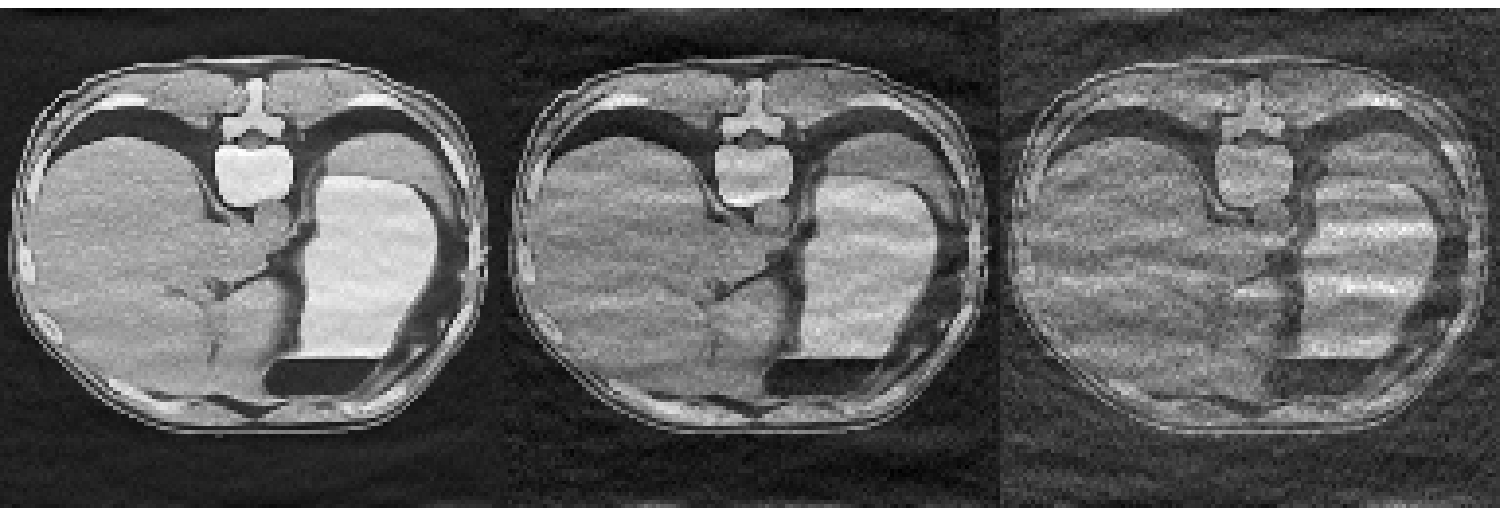}
  \caption{Low noise (left), moderate noise (middle), and higher noise (right).
  \label{noise}
 }
\end{figure}

In additional experiments we examined the effect of noise in our reconstruction.  The noise model we used is a simple stochastic model $|J|_{n} = |J| + \gamma*R$, where
$R$ is the normally distributed pseudo-random matrix of the order as $|J|$
with mean zero and standard deviation of one, and $\gamma > 0$ is the model
standard deviation chosen as $\gamma = \delta*|| |J| ||/||R||$,
where $\delta$ is the noise level, i.e., $(| J|_{n}- |J|)/|J|$. In Figure \ref{noise}, we show reconstructed images obtained by 20
iterations of the alternating split Bregman algorithm for three different levels of noise. The numerical values of the $l_2$-relative errors vs the exact solution for Figure \ref{noise} are shown in Talbel \ref{table3}.

\begin{table}[!h]
    \begin{tabular}{ | l | l | l | }
    \hline
      Low Noise (Level=0.01)   &  Moderate Noise ( Level=0.035) &  Higher Noise ( Level=0.06)  \\ \hline
      0.026       &     0.080      &    0.152 \\ \hline
    \end{tabular}
   \caption{Numerical error for alternating split Bregman}    
     \label{table3}    
   \end{table}

\subsection*{Conclusion}

We presented a convergent alternating split Bregman algorithm for least gradient problems with Dirichlet boundary data. This, in particular, leads to a convergent algorithm for recovering an isotropic conductivity $\sigma$ form the knowledge of the magnitude $|J|$ of one current generated by a given voltage $f$ on the boundary. Duality plays an essential role in our convergence proof and leads to a novel method to recover the full vector field $J$ from knowledge of its magnitude, and of the voltage $f$ on the boundary.  The alternating split Bregman algorithm presented here converges for non two-to-one boundary data as well as two-to-one boundary data. Several successful numerical experiments demonstrated the convergence behavior of the proposed algorithm.

 \end{document}